\documentclass{article}

\usepackage[german,english]{babel}
\usepackage{amsmath}
\usepackage{amsthm}
\usepackage{amssymb,latexsym}
\usepackage{url}
\usepackage{graphicx}
\usepackage{booktabs}
\usepackage{tabularx}
\usepackage{psfrag}

\newtheorem{deff}{Definition}
\newtheorem{lem}[deff]{Lemma}
\newtheorem{prop}[deff]{Proposition}
\newtheorem{thm}[deff]{Theorem}
\newtheorem{cor}[deff]{Corollary}

\newtheorem{prob}[deff]{Problem}

\makeatletter

\title{Knotted Polyhedral Tori}

\author{\Large Frank H.~Lutz\footnote{Supported by the DFG Research
     Group ``Polyhedral Surfaces'', Berlin, and by the MSRI-Program
     ``Computational Applications of Algebraic Topology'' (2006), Berkeley}\,\,
     and Nikolaus Witte$^*$}

\date{}

\begin{document}

\selectlanguage{english}

\maketitle

\begin{abstract}
  \noindent
  For every knot $K$ with stick number $k$ there is a knotted
  polyhedral torus of knot type $K$ with $3k$ vertices.
  We prove that at least $3k-2$ vertices are necessary.
\end{abstract}

\section{Introduction}

According to Alexander~\cite{Alexander1924}, every piecewise linear (PL) 
embedding of the $2$-torus $T^2$ in the $3$-sphere~$S^3$ splits $S^3$ into two
parts of which at least one part is a solid $3$-torus. The other part not
necessarily is a solid $3$-torus. As a consequence, a PL embedded $2$-torus 
in ${\mathbb R}^3$ need not bound a solid $3$-torus. If it does not, however, then
the one-point compactification of the non-compact component is a solid $3$-torus. 

Every PL embedded solid $3$-torus in ${\mathbb R}^3$ is isotopic to a tubular 
neighborhood of a knot $K$ in ${\mathbb R}^3$.
The knot type of $K$ is hereby fixed by the PL embedding of the solid $3$-torus.
The \emph{stick number} $s(K)$ of the knot $K$ is the minimal number of
\emph{sticks} (straight line segments) that are needed to built (an isotopic
copy of) $K$; see \cite{AdamsBrennanGreilsheimerWoo1997,Negami1991,Randell1994,Scharein1998}.

\begin{prob}\label{prob:1}
Given a knot $K$ with stick number $k$. How many vertices $P(K)$ are
at least needed for a polyhedral $2$-torus in ${\mathbb R}^3$ such that it bounds 
a solid $3$-torus isotopic to a tube around $K$?
\end{prob}

\noindent
We call such a polyhedral $2$-torus a \emph{knotted polyhedral torus
of knot type $K$}.

If a polyhedral $2$-torus in ${\mathbb R}^3$ bounds the compact complement of an 
once punctured solid $3$-torus of knot type $K$, then we call the $2$-torus 
a \emph{polyhedral torus of complement knot type~$K$}.
We say that it bounds a \emph{solid knot complement}.

\begin{prob}\label{prob:2}
Given a knot $K$ with stick number $k$. How many vertices $\overline{P}(K)$ are
at least needed for a polyhedral $2$-torus in ${\mathbb R}^3$ of complement knot type~$K$?
\end{prob}

A \emph{polyhedral map} on a (compact) surface $M$ is a decomposition 
of $M$ into an abstract polyhedral complex, i.e., a finite set of 
vertices, edges and polygons such that any two of the polygons 
either intersect in a common edge, a single vertex, or do not intersect at all 
(cf.\ \cite{BrehmSchulte1997,BrehmWills1993}).
In general, it is a hard problem (see \cite{BokowskiSturmfels1989})
to decide whether a given polyhedral map can be realized geometrically
as a polyhedron in $3$-space, i.e., with straight edges, convex faces, and without 
non-trivial intersections.

Every polyhedral map on $S^2$ is realizable as the boundary complex of a
convex $3$-polytope due to Steinitz \cite{Steinitz1922,SteinitzRademacher1934}. 
However, it is still an open problem, whether every triangulation of the $2$-torus 
is realizable in ${\mathbb R}^3$ (cf.\ Duke~\cite{Duke1970} and Gr\"unbaum~\cite[Ch.~13.2]{Gruenbaum1967}).
For recent work on the realizability of triangulated surfaces with few vertices 
and for further references see 
\cite{Bokowski2006pre,BokowskiGuedes_de_Oliveira2000,HougardyLutzZelke2006pre,Lutz2005pre,Schewe2007}.

A simulated annealing type heuristic was used in \cite{HougardyLutzZelke2006pre}
to obtain explicit realizations, in particular, of some examples of triangulated $2$-tori 
with up to $30$ vertices. All the resulting realizations of tori turned out to be 
unknotted, which arose our interest in finding a vertex-minimal example
of a knotted polyhedral torus.

Every $2$-face of a polyhedral torus with $n$ vertices in ${\mathbb R}^3$  
is a convex polygon. All the polygons with more than three sides
can be subdivided by adding diagonals to yield a geometrically realized 
triangulation of the torus with the same number $n$ of vertices. 
Any such resulting realization will have coplanar triangles. 
However, by perturbing the coordinates of the vertices slightly, 
we get a realization of the torus with vertices in general position 
(i.e., no three vertices on a line and no four vertices on a plane). 
Thus for the Problems~\ref{prob:1} and~\ref{prob:2} it suffices 
to consider geometric realizations (in general position) of triangulated tori.

Let us point out that there are no explicit tools available to
obtain knotted or linked realizations of triangulated surfaces.
The oriented matroid approach \cite{BokowskiSturmfels1989}
to realization problems allows (theoretically) to decide whether 
a given triangulated surface is realizable or not. 
However, due to its complexity, it is not applicable
in practice, and it is unclear how to suitably 
built in the knottedness requirement.

In Section~\ref{sec:upper} we discuss further preliminaries.
In particular, we present a (trivial) upper bound of $3k$ for $P(K)$ 
and an upper bound of $3k+4$ for $\overline{P}(K)$
(where, in the special case of the unknot, we have $P(\mbox{unknot})=\overline{P}(\mbox{unknot})=7$).
A lower bound of $3k-2$ for $P(K)$ and $\overline{P}(K)$ 
is proved in Section~\ref{sec:3_k} (for triangulations
that have an empty triangle) and in Section~\ref{sec:mk_torus}
(for triangulations without an empty triangle). Additional remarks 
are given in Section~\ref{sec:remarks}.

\section{Upper bounds and further preliminaries}
\label{sec:upper}

Given a knot $K$ with stick number $k$, we can easily
build a knotted polyhedral torus of knot type $K$
with $3k$ vertices: W.l.o.g.\ let $K$ be given as
a knotted polygon with $k$ edges (and $k$ vertices).
If we choose $\epsilon>0$ small enough, then the union 
of all $\epsilon$-balls with centers on $K$ is a tubular neighborhood of $K$.
Let $v$ be one of the $k$ vertices of $K$.
We replace $v$ by three new (distinct) vertices $v_1$, $v_2$, and~$v_3$
on the circle of radius $\epsilon$ that is the intersection 
of the boundary of the $\epsilon$-ball with center $v$ and the hyperplane 
that is the angle bisector at $v$. For every edge $v$--$w$ of $K$ 
consider the boundary of the convex hull of the six new vertices 
$v_1$, $v_2$, $v_3$, $w_1$, $w_2$, and $w_3$. 
If we remove from this boundary the two triangles $v_1v_2v_3$ and
$w_1w_2w_3$, we obtain a cylinder $C_{\text{$v$--$w$}}$. The union of
these cylinders $C_{\text{$v$--$w$}}$ for all edges $v$--$w$
of $K$ is a knotted polyhedral torus of knot type $K$ 
with $3k$ vertices.
\begin{prop}
Every knot $K$ with stick number $k$ can be modeled polyhedrally
by a knotted polyhedral torus of knot type $K$ with $3k$ vertices,
i.e., $P(K)\leq 3k$ for every knot $K$ with stick number $k$.
\end{prop}

For polyhedral tori of complement knot type~$K$
the following modification of the above construction
was pointed out to us by John M.\ Sullivan: 
W.l.o.g.\ let $K$ be given as a knotted polygon with $k$ edges 
and coordinates in general position.
Let further $v$ be a vertex of $K$ in convex position (i.e.,
$v$ is a vertex of the convex hull of the $k$ vertices of $K$).
We first proceed as above and build a knotted polyhedral torus
of knot type $K$ with $3k$ vertices. For this, the vertex $v$ 
is replaced by the three vertices $v_1$, $v_2$, and $v_3$.
Since $v$ was in convex position with respect to $K$,
w.l.o.g.\ we may assume that the edge $v_1$--$v_2$ is in convex position
with respect to the knotted polyhedral torus.
The edge $v_1$--$v_2$ lies in a triangle $v_1v_2w_i$, where $w_i$ is one of
the three vertices $w_1$, $w_2$, $w_3$ replacing a neighboring vertex $w$ of $v$. 
By placing a new vertex $y$ suitably ``above'' the triangle $v_1v_2w_i$
and close to the edge $v_1$--$v_2$ we obtain a knotted polyhedral torus
with subdivided triangle $v_1v_2w_i$ for which the triangle
$v_1v_2y$ is in convex position. Let now $z_1z_2z_3$ be a triangle
such that the convex hull of the six vertices $v_1$, $v_2$, $y$,
$z_1$, $z_2$, and $z_3$ (which is the boundary of an octahedron)
encloses the knotted polyhedral torus. If we remove from the union of the 
knotted polyhedral torus with the boundary of the octahedron
the triangle $v_1v_2y$, we obtain a polyhedral $2$-torus in ${\mathbb R}^3$ 
of complement knot type~$K$ with $3k+4$ vertices.

\begin{prop}[{John M.\ Sullivan}]
For a given knot $K$ with stick number $k$ at most $3k+4$ vertices are needed 
to built a polyhedral $2$-torus in ${\mathbb R}^3$ of complement knot type~$K$,
i.e., $\overline{P}(K)\leq 3k+4$ for every knot $K$ with stick number $k$.
\end{prop}

\begin{figure}
    \begin{center}
      \psfrag{1}{1}
      \psfrag{2}{2}
      \psfrag{3}{3}
      \psfrag{4}{4}
      \psfrag{5}{5}
      \psfrag{6}{6}
      \psfrag{7}{7}
      \includegraphics[height=33mm]{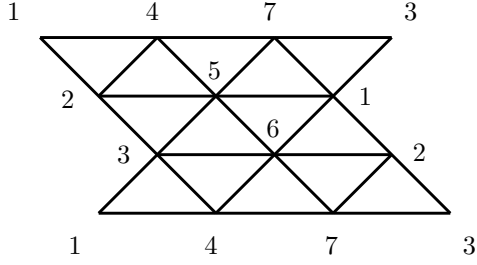}
    \end{center}
    \caption{M\"obius' $7$-vertex torus.
      \label{fig:moebius}}
  \end{figure}

The unique vertex-minimal triangulation of the $2$-torus is M\"obius' torus \cite{Moebius1886} 
with $7$ vertices; see Figure~\ref{fig:moebius}.
A polyhedral realization in ${\mathbb R}^3$ of M\"obius' torus was first
given by Cs\'asz\'ar \cite{Csaszar1949}. Bokowski and Eggert \cite{BokowskiEggert1991}
proved that there are essentially $72$ ``different'' polyhedral realizations of M\"obius' torus.
All of these bound an unknotted solid $3$-torus. Thus, $P(\mbox{unknot})=\overline{P}(\mbox{unknot})=7$. 
If $K$ is a knot (with stick number $k$) different from the unknot,
then $k\geq 6$ \cite{Randell1994} (with the trefoil knot
as the only knot of stick number exactly $6$).

Let $T$ be a triangulated torus and let $L$ be a cycle (i.e., a simple
closed curve) of $T$. We denote by $l_L$ the \emph{length} of $L$, that is,
the number of edges of $L$. A cycle $L$ of $T$ is \emph{non-separating}
if it does not bound a disc in $T$. Let a \emph{marked torus} 
be a pair $(T,M)$ consisting of a triangulated torus $T$ 
together with a non-separating cycle $M$ of $T$. For every choice of $M$ 
let $m_M$ be the length of a smallest cycle homotopic to $M$ and 
let $k_M$ be the length of a smallest non-separating cycle not homotopic to $M$. 
We call $m_M \times k_M$ the \emph{type} of the marked torus $(T,M)$.

\begin{deff}
For a given triangulated torus $T$ let the \emph{combinatorial stick number} $s(T)$ 
be the maximal value $k_M$ that is possible for some choice of a non-separating cycle $M$.
\end{deff}

Let $T$ be a triangulated torus with combinatorial stick number $s(T)$
and let $m$ be the length of a smallest non-separating cycle in $T$.
We call $m\times s(T)$ the \emph{type} of $T$.

\begin{lem}
Let $T$ be a triangulated torus of type $m\times s(T)$ and let
$M$ be a non-separating cycle of length $l_M=m$. Then $(T,M)$ is of type $m\times s(T)$.
\end{lem}

\begin{proof}
For a triangulated torus of type $m\times s(T)$ we have $m\leq s(T)$ 
(by definition of $m$ and $s(T)$).
If $M$ is a shortest non-separating cycle of length $l_M=m$,
then $s(T)$ is the length of a shortest non-separating cycle
not homotopic to $M$. Thus, there is at most one homotopy class
of non-separating cycles containing a cycle of length less than $s(T)$.
If $m<s(T)$, then this unique class contains $M$, and $(T,M)$ is of type $m\times s(T)$.
If $m=s(T)$, then $(T,M)$ is of type $s(T)\times s(T)$.
\end{proof}

Let $T$ be a triangulated torus, realized as a knotted polyhedral torus of knot type $K$. 
A non-separating cycle $M$ of $T$ is a \emph{meridian cycle} (in the realization of $T$) 
if $M$ is contractible in the solid $3$-torus bounded by $T$. Any two meridian cycles $M_1$ and $M_2$
are homotopic in $T$, and if $M'$ is a cycle of $T$ homotopic to a meridian cycle $M$,
then $M'$ is also a meridian cycle. (For polyhedral tori that bound a solid knot complement
meridians can be defined analogously.) 

\begin{lem}
If $N$ is a non-meridian non-separating cycle in a knotted polyhedral torus of knot type~$K$
(or in a polyhedral torus of complement knot type $K$)
with stick number $s(K)$, then $N$ has length $l_N\geq s(K)$.
\end{lem}

\begin{proof}
Any non-separating cycle $N$, which is not a meridian cycle, is isotopic to $K$
with stick number~$s(K)$ or isotopic to a satellite knot of $K$ with stick
number at least $s(K)$. Thus $l_N\geq s(K)$.
\end{proof}

\begin{cor}\label{cor:st_sk}
Let $T$ be a triangulated torus with combinatorial stick number $s(T)$.
If $T$ is realizable as a knotted polyhedral torus of knot type $K$
(or as a polyhedral torus of complement knot type $K$), then $s(T)\geq s(K)$.
\end{cor}
As a consequence, if $T$ is a knotted polyhedral torus of knot type
$K$ (or if $T$ is a polyhedral torus of complement knot type $K$) 
and if $T$ contains a non-separating cycle $N$ of length $l_N<s(K)$,
then necessarily $N$ is a meridian cycle.

\section{Tori of type \mathversion{bold}$3\times k$\mathversion{normal}}
\label{sec:3_k}

M\"obius' unique $7$-vertex torus is of type $3\times 3$
(as can easily be verified in Figure~\ref{fig:moebius})
and can be used to polyhedrally model the unknot with stick number $3$.

In the following, let $T$ be a triangulated torus with combinatorial stick number $s(T)$ 
and let $T$ have an \emph{empty triangle}, i.e., a non-separating cycle $M$ of length $l_M=3$.
Then $T$ is of type $3\times s(T)$.

\begin{thm}\label{thm:3k_torus}
  Any triangulated torus of type $3\times k$ has at least $3k-2$ vertices. 
  Furthermore, there is a unique vertex-minimal triangulation of type $3\times k$ 
  of the torus with $3k-2$ vertices.
\end{thm}

\begin{proof}
  We first show that any triangulated torus $T$ of type $3\times k$ has at least $3k-2$
  vertices. If we cut open $T$ along one of the non-separating cycles of length~$3$, 
  we obtain a cylinder, as depicted in bold in Figure~\ref{fig:cylinder}.
  \begin{figure}
    \begin{center}
      \psfrag{1}{1}
      \psfrag{2}{2}
      \psfrag{3}{3}
      \psfrag{I}{I}
      \psfrag{II}{II}
      \psfrag{a}{$b$}
      \psfrag{b}{$c$}
      \psfrag{c}{$a$}
      \includegraphics[height=33mm]{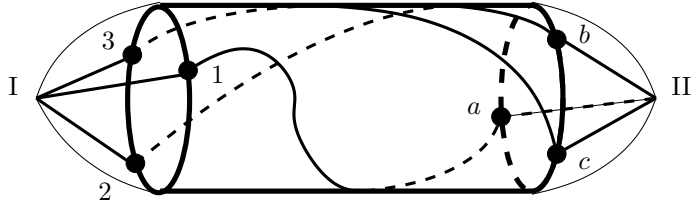}
    \end{center}
    \caption{The cylinder and two cones.
      \label{fig:cylinder}}
  \end{figure}
  Let the vertices on the ``left'' side of the cylinder be $1$, $2$, and $3$, 
  and the vertices on the ``right'' side of the cylinder be $a$, $b$, and $c$, 
  in circular order respectively. We close the cylinder by adding two cones, 
  one with apex I over the circle 1--2--3--1 and one with apex II 
  over the circle \mbox{$a$--$b$--$c$--$a$}, to obtain a triangulated $2$-sphere. 
  By Steinitz' theorem \cite{Steinitz1922,SteinitzRademacher1934}, 
  every triangulated $2$-sphere is realizable as the boundary complex 
  of a convex $3$-polytope. Moreover, the $1$-skeleton of a $3$-polytope 
  is a $3$-con\-nected (planar) graph, i.e., for every pair of its vertices, 
  there are three independent paths in the graph connecting the two vertices. 
  For every possible choice of three independent paths that connect the vertices I and II, 
  the respective paths must use the edges I--1, I--2, I--3, 
  $a$--II, $b$--II, and $c$--II, but cannot use one of the edges 
  1--2, 2--3, 3--1, $a$--$b$, $b$--$c$, and $c$--$a$.  
  W.l.o.g.\ let the paths connect vertex 1 with $a$, 2 with
  $b$, and $3$ with $c$, respectively, as in Figure~\ref{fig:cylinder}.  
  If we cut open the cylinder along the path $1$---$a$ (see Figure~\ref{fig:threepaths}), 
  then there are three ways to glue back together the torus $T$ by identifying the vertices
  $1$, $2$, and $3$ on the left hand side with the vertices $a$, $b$, and $c$
  on the right hand side:
  \begin{itemize}
  \item[a)] $1':=a\equiv 1$, $2':=b\equiv 2$, $3':=c\equiv 3$,
  \item[b)] $3':=a\equiv 3$, $1':=b\equiv 1$, $2':=c\equiv 2$,
  \item[c)] $2':=a\equiv 2$, $3':=b\equiv 3$, $1':=c\equiv 1$.
  \end{itemize}

  \begin{figure}
    \begin{center}
      \psfrag{1}{1}
      \psfrag{2}{2}
      \psfrag{3}{3}
      \psfrag{4}{4}
      \psfrag{5}{5}
      \psfrag{6}{6}
      \psfrag{7}{7}
      \psfrag{8}{8}
      \psfrag{9}{9}
      \psfrag{10}{10}
      \psfrag{11}{11}
      \psfrag{12}{12}
      \psfrag{x}{$3k-5$}
      \psfrag{y}{$3k-4$}
      \psfrag{z}{$3k-3$}
      \psfrag{a}{$a=3'$}
      \psfrag{b}{$b=1'$}
      \psfrag{c}{$c=2'$}
      \includegraphics[height=33mm]{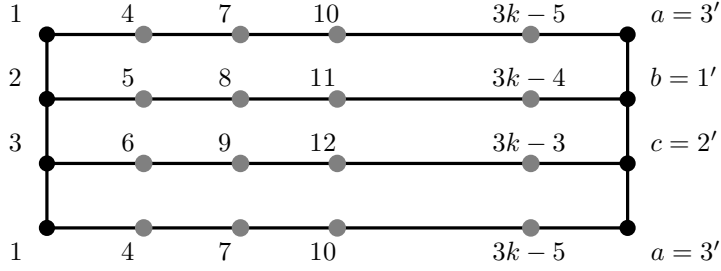}
    \end{center}
    \caption{The cut open torus with three independent paths.
      \label{fig:threepaths}}
  \end{figure}

  In order to complete Figure~\ref{fig:threepaths} to a triangulation
  of type $3\times k$ of the torus, every connecting path 1---1',
  2---2', and 3---3' must have at least length $k$.  Therefore, in
  the case a), each of the paths 1---$a$, 2---$b$, and 3---$c$ must be
  subdivided at least $k-1$ times, which gives together at least $3k$
  vertices for a respective triangulation of the torus, so we do not
  consider this case further.

  By symmetry, the case c) is equivalent to the case b), so let us analyze case b).
  We assume that the paths 1---$a$, 2---$b$, and 3---$c$ each are subdivided 
  the minimal number of times $k-2$, as indicated in Figure~\ref{fig:threepaths},
  which yields in total at least $3k-3$ vertices for a respective triangulation 
  of the torus. If we further assume that there are no additional vertices, i.e., 
  there are no interior vertices within the three polygonal strips, then we will see
  in the following that the combinatorial stick number of~$T$ is less than $k$, a contradiction.
  Thus, at least one more vertex is needed, which will give the 
  lower bound of $3k-2$ vertices.

  In every triangulation of a polygon $P$ without interior vertices
  there are at least two vertices that are not the endpoints of a
  diagonal.  If for three consecutive vertices $j-1$, $j$, and $j+1$
  of $P$ the vertex~$j$ is not the endpoint of a diagonal, then the
  edge $(j-1)$--$(j+1)$ is present in the triangulation.

  In case b) we are obviously not allowed to add edges of type
  $i$--$(i+6)$ to Figure~\ref{fig:threepaths}, since this would yield
  a path 1---1', 2---2', or 3---3' of length less then $k$. Also we
  are not allowed to include the edges 1--5, 2--6, 3--4, ($3k-5$)--1',
  ($3k-4$)--2', and $(3k-3)$--3'. This leaves the vertices 1 and 1' in
  the upper strip, 2 and 2' in the middle strip, and 3 and 3' in the
  lower strip as vertices that are not endpoints of a diagonal. In
  other words, the edges 2--4, $(3k-4)$--3', 3--5, $(3k-3)$--1', 1--6,
  and $(3k-5)$--2' have to be present in any resulting triangulation
  and therefore reduce the polygonal strips to strips that have two
  vertices less, respectively. If we apply the same line of arguments
  to the strips with two vertices less, we are forced to uniquely add 
  further edges. We obtain a unique resulting triangulation as depicted
  in Figure~\ref{fig:threepaths_b}.  Unfortunately, in this triangulation
  there are paths connecting 1---1' (and also 2---2', 3---3') of
  length $k-1$, e.g., the path 1--6--8--11--14--\dots--$(3k-4)$--1'.
  \begin{figure}
    \begin{center}
      \psfrag{1}{1}
      \psfrag{2}{2}
      \psfrag{3}{3}
      \psfrag{4}{4}
      \psfrag{5}{5}
      \psfrag{6}{6}
      \psfrag{7}{7}
      \psfrag{8}{8}
      \psfrag{9}{9}
      \psfrag{10}{10}
      \psfrag{11}{11}
      \psfrag{12}{12}
      \psfrag{x}{$3k-5$}
      \psfrag{y}{$3k-4$}
      \psfrag{z}{$3k-3$}
      \psfrag{a}{3'}
      \psfrag{b}{1'}
      \psfrag{c}{2'}
      \includegraphics[height=33mm]{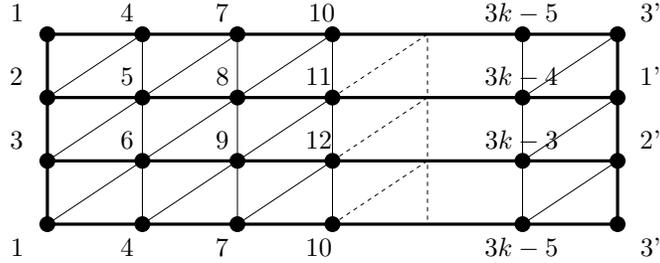}
    \end{center}
    \caption{The cut open torus of case b) with $3k-3$ vertices.
      \label{fig:threepaths_b}}
  \end{figure}
  It follows that there is no triangulation
  of a torus of type $3\times k$ with
  $3k-3$ (or fewer) vertices.

  We finally show that, up to isomorphy, there is exactly one triangulated torus of type
  $3\times k$ with $3k-2$ vertices. Such a triangulation, if it exists, can only be obtained 
  from case b) (or, by symmetry, from case c)) of Figure~\ref{fig:threepaths} 
  either 
  \begin{itemize}
  \item[b')] by subdividing once one of the paths, say, 1---$3'$, or 
  \item[b'')] by allowing an interior vertex in one of the polygonal
    strips, say, in the upper strip.
  \end{itemize}

  In the case b'') the middle and the lower strip are triangulated as before. 
  If the additional vertex~$x$ in the upper strip is only connected to some 
  of the upper vertices $1,4,7,\dots,3k-5,3'$, then this obviously forces a shortcut 
  for the path $1$---$3'$--$1'$ of length $k$. Similarly, if $x$ is only connected 
  to some of the lower vertices $2,5,8,\dots,3k-4,1'$ of the upper polygon,
  then we get a shortcut for the path $2$---$1'$--$2'$. Thus $x$ is connected
  to at least one upper and to at least one lower vertex.
  These vertices cannot lie too far apart. For example,
  if the additional vertex $x$ is connected to $5$ and $10$, 
  then the dotted path $5$--$x$--$10$ yields a shortcut of length $k-1$
  for the path $3$---$3'$. Hence, the only
  admissible cases are $1$--$x$--$5$, $4$--$x$--$8,\dots,(3k-5)$--$x$--$1'$ 
  as well as $2$--$x$--$4$, $5$--$x$--$7,\dots,(3k-4)$--$x$--$3'$.
  By the same arguments as above, however, then at least one of the
  vertices 1 or $1'$ is not the endpoint of a diagonal.
  In other words, at least one of the edges $2$--$4$ or
  $(3k-4)$--$3'$ is present in a resulting triangulation,
  each of which yielding a shortcut.

  \begin{figure}
    \begin{center}
      \psfrag{1}{1}
      \psfrag{2}{2}
      \psfrag{3}{3}
      \psfrag{4}{4}
      \psfrag{5}{5}
      \psfrag{6}{6}
      \psfrag{7}{7}
      \psfrag{8}{8}
      \psfrag{9}{9}
      \psfrag{10}{10}
      \psfrag{11}{11}
      \psfrag{12}{12}
      \psfrag{x}{$3k-5$}
      \psfrag{y}{$3k-4$}
      \psfrag{z}{$3k-3$}
      \psfrag{a}{3'}
      \psfrag{b}{1'}
      \psfrag{c}{2'}
      \psfrag{s}{$x$}
      \includegraphics[height=33mm]{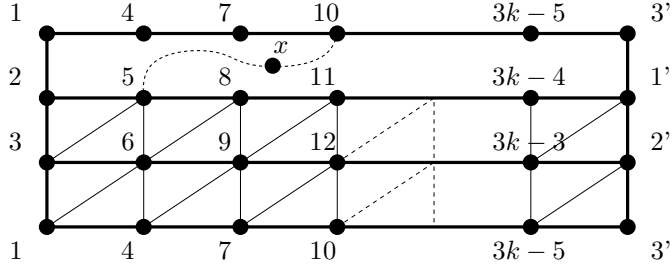}
    \end{center}
    \caption{The cut open torus of case b'') with interior vertex.
      \label{fig:threepaths_interior}}
  \end{figure}

  Thus we are left with case b') for which there is
  a unique resulting triangulation as displayed in 
  Figure~\ref{fig:3k_2_torus} (with the middle strip as before
  and unique ways to triangulate the upper and lower strips).

  \begin{figure}
    \begin{center}
      \psfrag{1}{1}
      \psfrag{2}{2}
      \psfrag{3}{3}
      \psfrag{4}{4}
      \psfrag{5}{5}
      \psfrag{6}{6}
      \psfrag{7}{7}
      \psfrag{8}{8}
      \psfrag{9}{9}
      \psfrag{10}{10}
      \psfrag{11}{11}
      \psfrag{12}{12}
      \psfrag{3k-5}{$3k-5$}
      \psfrag{3k-4}{$3k-4$}
      \psfrag{3k-3}{$3k-3$}
      \psfrag{3k-2}{$3k-2$}
      \psfrag{a}{3}
      \psfrag{b}{1}
      \psfrag{c}{2}
      \includegraphics[height=33mm]{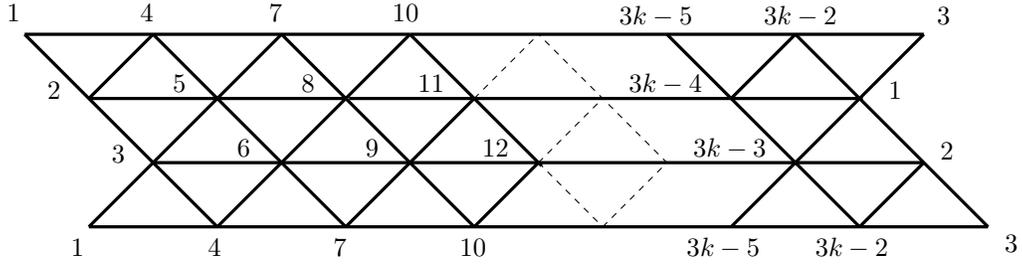}
    \end{center}
    \caption{The unique triangulated torus of type $3\times k$ with $3k-2$ vertices.
      \label{fig:3k_2_torus}}
  \end{figure}
\end{proof}

\section{Tori of type \mathversion{bold}$m\times k$\mathversion{normal}}
\label{sec:mk_torus}

In this section, we give a lower bound for the number of vertices of
triangulated tori of general type $m\times k$. For $k\geq 6$ this
bound is at least $3k-2$. As a consequence, for every knot $K$ with
stick number~$k$, $P(K)$ and $\overline{P}(K)$ are bounded from below
by $3k-2$.

Let $T$ be a triangulated torus of type $m\times k$ and let~$V$ be the
set of vertices of $T$. Let the distance ${\rm dist}(v,w)$ of two
vertices $v,w\in V$ be the length of a shortest path connecting the
vertices $v$ and $w$.  For $i\geq 0$, let $V_i(v)=\{w\in V\,|\,{\rm dist}(v,w)=i\}$ 
be the set of vertices of the torus that have distance $i$ from $v$.  
Since $T$ has combinatorial stick number~$k$,
$V_i\neq\emptyset$ for $0\leq i\leq\left\lfloor\frac{k}{2}\right\rfloor$. 

\begin{prop}\label{prop:mk_torus}
  Let $T$ be a triangulated torus of type $m\times k$.  Then
  \[|V|\geq 2\left\lceil\frac{m}{2}\right\rceil^2 +
  \big(k-2\left\lceil\frac{m}{2}\right\rceil\big)m + 1.\]
\end{prop}

\begin{proof}
  Let $T$ be a triangulated torus of type $m\times k$, let $M$ be a
  (minimal) non-separating cycle of $T$ of length $m$, and let $v$ be a vertex
  on $M$. We bound the number of vertices of each~$V_i$ for $0\leq
  i\leq\left\lfloor \frac{k}{2} \right\rfloor$ from below.

  For $1\leq i\leq\left\lfloor \frac{k-1}{2} \right\rfloor$ each~$V_i$ is split 
  by~$M$ into two sets which we call the ``left'' and the ``right''
  part of~$V_i$; see Figure~\ref{fig:mk_torus}. The vertices of~$V_i$,
  which lie on~$M$, contribute to the ``right'' part of~$V_i$. 
  Let~$A_0=V_0$. For~$1\leq i\leq \left\lceil \frac{m}{2} \right\rceil -1$ we denote the ``right''
  part of~$V_i$ by~$A_i$ and for~$1\leq i\leq \left\lceil \frac{m}{2}
  \right\rceil$ the ``left'' part by~$E_i$. Similarly, we denote 
  for~$\left\lceil \frac{m}{2} \right\rceil \leq i\leq \left\lfloor \frac{k-1}{2} \right\rfloor$ 
  the ``right'' part of~$V_i$ by~$B_i$ and 
  for~$\left\lceil \frac{m}{2} \right\rceil+1 \leq i\leq \left\lfloor \frac{k-1}{2} \right\rfloor$
  the ``left'' part of~$V_i$ by~$D_i$
. In the case that~$m$ is even, we set~$C_{k/2}=V_{k/2}$.
    \begin{figure}
    \begin{center}
      \psfrag{v}{$v$}
      \psfrag{A0}{$A_0$}
      \psfrag{A1}{$A_1$}
      \psfrag{A2}{$A_2$}
      \psfrag{A3}{$A_3$}
      \psfrag{B4}{$B_4$}
      \psfrag{B5}{$B_5$}
      \psfrag{C6}{$C_6$}
      \psfrag{D5}{$D_5$}
      \psfrag{E4}{$E_4$}
      \psfrag{E3}{$E_3$}
      \psfrag{E2}{$E_2$}
      \psfrag{E1}{$E_1$}
      \psfrag{M}{$M$}
      \includegraphics[height=70mm]{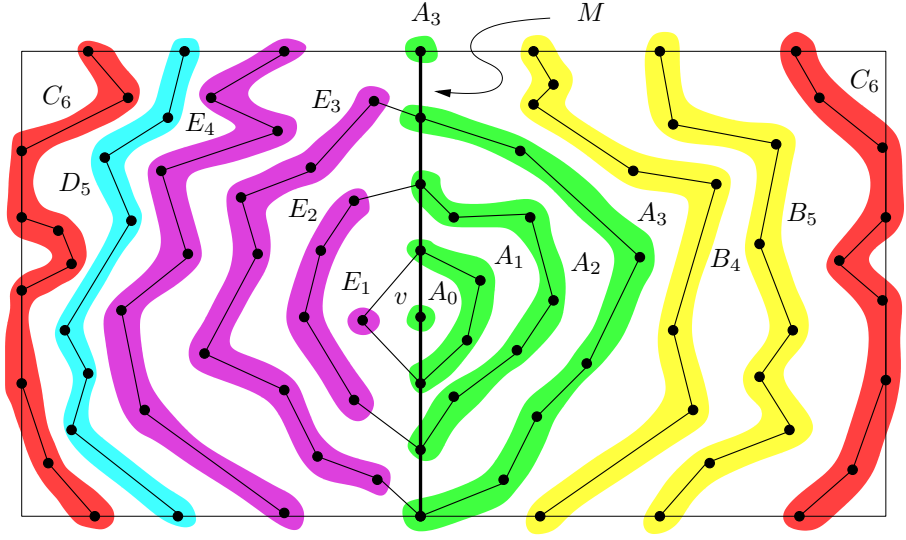}
    \end{center}
    \caption{A torus of type $7\times 12$.
      \label{fig:mk_torus}}
  \end{figure}
  It follows that
  \begin{equation}\label{eq:mk_torus0}
    \begin{aligned}
      |A_i| \geq 2i+1  \quad &\text{for} \quad  0\leq i\leq \lceil m/2 \rceil -1,\\
      |B_i| \geq m     \quad &\text{for} \quad  \lceil m/2 \rceil \leq i\leq \lfloor (k-1)/2 \rfloor,\\
      |C_{k/2}| \geq m \quad &\text{in the case~$m$ even,}\\
      |D_i| \geq m     \quad &\text{for} \quad  \lceil m/2 \rceil +1\leq i\leq \lfloor (k-1)/2 \rfloor\text{, and}\\
      |E_i| \geq 2i-1  \quad &\text{for} \quad  1\leq i\leq \lceil m/2 \rceil,
    \end{aligned}
  \end{equation}
  since violating any of the inequalities above would imply the
  existence of a non-separating cycle of length less than~$m$.
  The lower bound stated in Proposition~\ref{prop:mk_torus} is
  obtained by summing over the lower bounds for the number of vertices
  of the~$A_i$'s,~$B_i$'s,~$C_{k/2}$ (for~$m$ even),~$D_i$'s, and~$E_i$'s. There
  are~$\left\lceil\frac{m}{2}\right\rceil$ summands each for the~$A_i$'s and
  the~$E_i$'s, and~$k-2\left\lceil\frac{m}{2}\right\rceil$
  remaining summands corresponding to the~$B_i$'s,~$C_{k/2}$ (for~$m$ even),
  and~$D_i$'s. Overall we accumulate
  \begin{equation}\label{eq:mk_torus}
    \begin{aligned}
      |V| &\geq
      2 \sum_{j=1}^{\left\lceil\frac{m}{2}\right\rceil}(2j+1)
      \;+\;
      \sum_{j=1}^{k-2\left\lceil\frac{m}{2}\right\rceil} m\\
      &=2\left\lceil\frac{m}{2}\right\rceil^2
      \;+\;
      \big(k-2\left\lceil\frac{m}{2}\right\rceil\big)m.
    \end{aligned}
  \end{equation}

  Finally, we prove that~$T$ must have at least one more vertex. To this end,
  assume that Inequality~\eqref{eq:mk_torus} is tight for~$T$. Then each
  inequality in~\eqref{eq:mk_torus0} must be tight. In particular,
  we have~$|A_1|+|E_1|=4$ and therefore~$v$ has degree four. Now~$v$ may be
  chosen arbitrarily on~$M$, thus all vertices on~$M$ have degree four.
  Let~$w_1$ be one of the two neighboring vertices of~$v$ on~$M$ 
  and let~$w_2\neq v$ be the other neighbor of~$w_1$ on~$M$. Further, let~$w_3$ be 
  the unique vertex of~$A_1$ not in~$M$. Since~$w_2$--$w_3$ cannot
  be an edge of~$T$ (it would create a non-separating cycle of length~$m-1$),
  $w_1$ must have degree at least five. Contradiction.
\end{proof}

\begin{thm}\label{thm:mk_torus}
  For a knot $K$ with stick-number~$k$, $P(K)$ and~$\overline{P}(K)$ 
  are bounded from below by~$3k-2$. Moreover, there is a unique
  triangulated torus on~$3k-2$ vertices with combinatorial stick-number~$k$.
\end{thm}

\begin{proof}
  By Corollary~\ref{cor:st_sk} the stick number $s(T)$ of a realized triangulated torus $T$ 
  is bounded from below by $k=s(K)$. Theorem~\ref{thm:3k_torus} proves 
  the existence of a unique vertex-minimal triangulated torus 
  of type~$3\times k$ on~$3k-2$ vertices.
  As mentioned in Section~\ref{sec:upper}, we have
  $P(\mbox{unknot})=\overline{P}(\mbox{unknot})=7=3\cdot 3-2$ for the unknot
  with stick number $3$. All other knots have stick number $k\geq 6$ \cite{Randell1994}.
  We show that for~$k\geq6$ and $m\geq4$ the lower
  bound in Proposition~\ref{prop:mk_torus} exceeds~$3k-2$, that is,
  \begin{equation*}
  2\left\lceil \frac{m}{2} \right\rceil^2 + (k-2\left\lceil\frac{m}{2}\right\rceil)m + 1\,\,  >\,\,  3k - 2
  \end{equation*}
  or, equivalently,
  \begin{equation}\label{eq:mk_torus2}
      (m-3)k + 2 \left\lceil \frac{m}{2} \right\rceil^2 - 2\left\lceil \frac{m}{2}
      \right\rceil m + 3\,\,  >\,\,  0,
  \end{equation}
which, for fixed~$m\geq4$, is a linear inequality in~$k$ with a
positive coefficient $(m-3)$ for~$k$. 
It is easy to verify that~\eqref{eq:mk_torus2} holds for~$m=4,5$
and~$k=6$, and for~$m=k\geq6$, thus completing the proof.
\end{proof}

\section{Remarks}
\label{sec:remarks}

As mentioned in the introduction, it is not known whether every triangulated torus 
can be realized geometrically in ${\mathbb R}^3$. However, the series of triangulated tori 
of type $3\times k$ with $3k-2$ vertices from Section~\ref{sec:3_k} \emph{is} realizable: 
The examples of the series have vertex-transitive cyclic symmetry, generated by the cycle
$(1,4,7,10,\dots,(3k-2),3,6,9,\dots,(3k-3),2,5,8,11,(3k-4))$,
and for $k\geq 3$, the example with $3k-2$ vertices is realizable in the boundary complex 
of the cyclic $4$-polytopes $C_4(3k-2)$; see \cite{Altshuler1971,KuehnelLassmann1985-di,KuehnelLassmann1996-bundle}.
(The examples are equivelar triangulations, i.e., all vertices have the same degree $6$, 
and equivelar triangulations of the torus all are vertex-transitive; see \cite{DattaUpadhyay2005}
and also \cite{BrehmKuehnel2006pre}.) The realization in $C_4(3k-2)$ is unknotted.

In generalization of our discussion of knotted realizations, one may ask for realizations 
in a given isotopy class. For example, if we cut open a torus,
realized as the boundary of a solid $3$-torus in ${\mathbb R}^3$, along a meridian, 
twist one end of the solid torus $t$ times and glue back together both ends, is the 
resulting embedded $2$-torus then realizable in the respective isotopy class? 

In the vertex-minimal torus of type $3\times k$ from Section~\ref{sec:upper}
the non-separating cycle
$1$--$4$--$7$--$10$--$\dots$--$(3k-2)$--$3$--$6$--$9$--$\dots$--$(3k-3)$--$2$--$5$--$8$--$11$--$(3k-4)$--$1$
is a Hamiltonian cycle. If the torus were realizable as a knotted polyhedral torus of knot type $K$ 
with stick number $k$, then the induced realization of the Hamiltonian cycle with $3k-2$ vertices 
would give a realization of a satellite knot of $K$. Is it possible to build this satellite knot 
with $3k-2$ sticks?

\subsection*{Acknowledgment}

The authors are grateful to J\"urgen Bokowski
and John M.\ Sullivan for helpful discussions
and comments.

\bibliography{.}

\bigskip
\medskip

\noindent
Frank H. Lutz\\
Technische Universit\"at Berlin\\
Institut f\"ur Mathematik, Sekr.\ MA 3-2\\
Stra\ss e des 17.\ Juni 136\\
10623 Berlin\\
Germany\\
{\tt lutz@math.tu-berlin.de}

\bigskip
\medskip

\noindent
Nikolaus Witte\\
Technische Universit\"at Berlin\\
Institut f\"ur Mathematik, Sekr.\ MA 6-2\\
Stra\ss e des 17.\ Juni 136\\
10623 Berlin\\
Germany\\
{\tt witte@math.tu-berlin.de}

\end{document}